\documentclass[psamsfonts]{amsart}
\usepackage{amssymb,latexsym,eucal,amsthm,amscd}
\usepackage[usenames]{color}

\usepackage[cp1250]{inputenc} 
\usepackage[T1]{fontenc}      

\definecolor{Black}{RGB}{0,0,0}
\definecolor{Red}{RGB}{255,0,0}
\definecolor{Blue}{RGB}{0,0,255}
\definecolor{Green}{RGB}{0,255,0}
\definecolor{Gray}{RGB}{120,120,120}

\theoremstyle{plain}                                     
\newtheorem{lmm}{Lemma}[section]
\newtheorem{thm}[lmm]{Theorem}
\newtheorem{prp}[lmm]{Proposition}
\newtheorem{crl}[lmm]{Corollary}

\newtheorem{clm}{Claim}

\theoremstyle{plain}
\newtheorem{prb}{Problem}[section]

\theoremstyle{definition}

\newenvironment{prf}
{\begin{proof}[Proof]\setcounter{clm}{0}}
{\end{proof}}

\newenvironment{cprf}
{\begin{proof}[Proof of Claim~\theclm]}
{ \end{proof}}
   
\newcommand{\Lmm}[1]{Lemma~\ref{#1}}
\newcommand{\Prp}[1]{Proposition~\ref{#1}}

\newcommand{\Crl}[1]{Corollary~\ref{#1}}

\newcommand{\Math}[1]{{\ensuremath{#1}}}  
\newcommand{\mbb}[1]{\Math{\mathbb{#1}}}

\newcommand{\ga}{\Math{\alpha}}
\newcommand{\gb}{\Math{\beta}}

\newcommand{\gl}{\Math{\lambda}}

\newcommand{\go}{\Math{\omega}}

\newcommand{\euler}[1]{\Math{\mathcal{#1}}}

\newcommand{\eP}{\euler{P}}

\newcommand{\NN}{\mbb{N}}

\DeclareMathOperator*{\Span}{Span}
\DeclareMathOperator*{\Spec}{Spec}
\DeclareMathOperator*{\Ann}{Ann}

\numberwithin{equation}{section}

\begin{document}

\title[]{Regularly weakly based modules over right perfect rings and Dedekind domains}
 
\author{Michal Hrbek, Pavel R\r{u}\v{z}i\v{c}ka}

\address{Charles University \\
          Faculty of Mathematics and Physics\\
          Department of Algebra \\
          Sokolovsk\'a 83\\
          186 75 Prague\\
          Czech Republic}     
          
\email{hrbmich@gmail.com, ruzicka@karlin.mff.cuni.cz}
\date{\the \day/\the \month/\the \year}

\thanks{ The first author is partially supported by the project SVV-2015-260227 of Charles University in Prague. Both the authors were partially supported by the Grant Agency of the Czech Republic under the grant no. GACR 14-15479S.} 

\newcommand{\Label}[1]{\label{#1}}

\begin{abstract} A \emph{weak basis} of a module is a generating set of the module minimal with respect to inclusion. A module is said to be \emph{regularly weakly based} provided that each of its generating sets contain a weak basis. In the paper we study  
\begin{enumerate}
\item rings over which all modules are regularly weakly based, refining results of Nashier and Nichols,
\item regularly weakly based modules over Dedekind domains. 
\end{enumerate} 
\end{abstract}

\maketitle

\section{Introduction}

By a module we always mean a right module over a ring $R$. Let $M$ a module and let $X, Y$ by subsets of $M$. We say that the set $X$ is \emph{weakly independent over $Y$} if $x \not\in \Span((X \setminus \{x\}) \cup Y)$ for all $x \in X$. We say shortly that $X$ is \emph{weakly independent} in the case of $Y=\emptyset$. A generating weakly independent subset of a module $M$ is called a \emph{weak basis} of $M$. A module $M$ is \emph{weakly based} if it contains a weak basis. Finally, a module $M$ is called \emph{regularly weakly based} if any generating set of $M$ contains a weak basis.

Nashier and Nichols characterized right perfect rings as rings over which every quasi-cyclic right $R$-module (i.e. every finitely generated submodule is contained in a cyclic one) is cyclic (i.e. every submodule is contained in a cyclic one). As a consequence of this they have got that rings over which all right modules are regularly weakly based are necessarily right perfect and they raised a question whether, conversely, all modules over right prefect rings are regularly weakly based. We refine their result proving that infinitely generated free modules over non-right perfect rings are not regularly weakly based and we observe that their question regarding right perfect rings easily reduces to semisimple rings. 

The other topic of the paper is a study of regularly weakly based modules over Dedekind domains. This is motivated by the characterization of weakly based modules over, first abelian groups \cite{HR1} and then Dedekind domains \cite{HR2} done by the authors. For regularly weakly based modules we will not obtain 
the full characterization, however we reduce the problem to a question of characterizing regularly weakly based modules over commutative semisimple rings, which indeed is a special case of the more general open question regarding right perfect rings introduced above.      
	
There is a few simple facts regarding regularly weakly based modules which we will freely use within the paper. Namely, it is clear that a finitely generated module is regularly weakly based. Also observe that unlike in case of weakly based modules, a direct summand of a regularly weakly based module is regularly weakly based. Also the next elementary lemma, in different variations, will be repeatedly used with no reference. Its proof is left to the reader.  

\begin{lmm}\label{L00c} Let $R$ be ring and $M$ a right $R$-module. Let $X,Y,Z$ be subsets of $M$. Suppose that $X$ is weakly independent over $Y \cup Z$ and
$Y$ is weakly independent over $X \cup Z$. Then $X \cup Y$ is weakly independent over $Z$.  
\end{lmm}

\section{Modules over right perfect rings}\label{Perfect}

We start with a natural task of characterizing rings $R$ such that all right $R$-modules are regularly weakly based. We refine the result of \cite{NN} that all such a modules must be right perfect. In particular, in \Lmm{L04}, we prove that an infinitely generate free module over a non-perfect ring is not regularly weakly based. Nashier and Nichols suggested conversely, that all modules over right perfect rings are regularly weakly based. We discuss this question in the final part of this section, adding an observation that we can factor out the Jacobson radical, and so reduce the question to semisimple rings.    

\begin{lmm}{\cite[Proposition 1 and Theorem 2]{NN}} \label{L03}
A ring $R$ is right perfect if and only if for each sequence $(r_n \mid n \in \omega)$ of elements of $R$ there is $n_0 \in \omega$ such that for all $n \geq n_0$ there is $j \geq 1$ such that $r_{n+j} \cdots r_{n+1} R = r_{n+j} \cdots r_{n+1} r_n R$.
\end{lmm}

\begin{lmm} \label{L04} Let $R$ be a ring that is not right perfect. Than a free right $R$-module is regularly weakly based if and only if it is finitely generated.
\end{lmm}
\begin{prf} A finitely generated module is regularly weakly based. Thus it suffices to show that an infinitely generated right free module is not regularly weakly based. Since a direct summand of a regularly weakly based module is regularly weakly based, we can restrict ourselves to a countably generated free right $R$ module $F=R^{(\aleph_0)}$. In order to show that $F$ is not regularly weakly based, fix a free basis $\{b_n \in n \in \omega\}$ of $F$. Since $R$ is not right perfect, there is by \Lmm{L03} a sequence $(r_n \mid n \in \omega)$ of elements of $R$ such that for any $n \in \omega$ and all $j \ge 1$ we have that $r_{n+j} \cdots r_{n+1} R \supsetneq r_n \cdots r_{n+1} r_{n} R$. In particular, this implies that all $r_n$ are non-invertible in $R$. For each $n \in \omega$, we define the following elements from $F$:
\begin{equation*} x_n = b_{n+1} r_n  \quad \text{and}\quad y_n = b_n - x_n = b_n - b_{n+1} r_n.
\end{equation*}	
Put $Z=\{x_n, y_n \mid n \in \omega\}$ and $Y = \{y_n \mid n \in \omega\}$. Clearly $B \subseteq \Span(Z)$, hence $Z$ generates $F$. We claim that $Z$ does not contain any weak basis of $F$. Suppose otherwise and pick a weak basis $W \subseteq Z$ of $F$. As $r_n$ are non-invertible, $Y \subseteq W$. 

Let $n \in \go$, and suppose $x_n \in W$. Observe that then $b_k$ and thus also $x_k = b_k - y_k$ belong to $\Span(W)$ for all $k \le n$. Since $W$ is weakly independent and $Y \subseteq W$, it contains at most one $x_n$, that is, $W \subseteq Y \cup \{x_n\}$ for some $n \in \go$. We claim that $b_{n+1} \notin \Span(W)$. Indeed, otherwise
\begin{equation}\label{Eq01} 
b_{n+1} = x_n s + \sum_{i \in \go} y_i s_i
\end{equation}
for some $s,s_0,s_1 \dots$ from $R$ such that all but finitely many $s_i$ are $0$. Substituting for $x_n$, $y_i$ we get that  	
\begin{equation}\label{Eq02} 
b_{n+1} = b_{n+1} r_n s + \sum_{i \in \go} (b_is_i - b_{i+1}r_is_i).
\end{equation}	
From this we get that $s_0 = \cdots = s_n = 0$, $s_{n+1} = (1 - r_n s)$ and $s_{n+1+j} = r_{n+j} \cdots r_{n+1} (1 - r_n s)$, for all $j > 0$. Since all but finitely many $s_i$ equal $0$, there is $j > 0$ such that $s_{n+1+j} = 0$. Then we get that $r_{n+j} \cdots r_{n+1} = r_{n+j} \cdots r_{n+1} r_n s$, which gives $r_{n+j} \cdots r_{n+1} R = r_{n+j} \cdots r_{n+1} r_n R$. This contradicts our choice of the sequence $(r_n \mid n \in \omega)$.   	
\end{prf}

Recall that a subset $I$ of a ring $R$ is right T-nilpotent provided that for every sequence $a_1,a_2,\dots$ there is a positive integer $n$ such that
$a_n \cdots a_1 = 0$. A right ideal $J$ of a ring $R$ is $T$-nilpotent if and only if $MJ \neq M$ for every non-zero right $R$-module \cite[Lemma 28.3]{AF}.   
By the Theorem of Bass \cite[Theorem 28.4]{AF} a ring $R$ is right perfect if and only if its Jacobson radical $J$ is $T$-nilpotent and the ring $R / J$ is semisimple.   

\begin{lmm}\label{L09} Let $J$ be a right $T$-nilpotent right ideal of a ring $R$, let $M$ be a right $R$-module. Then every $X \subseteq M$ lifting a weak basis of $M / M J$ over $MJ$ is a weak basis of $M$.
\end{lmm}

\begin{prf} Since $X$ lifts a weak basis of $M / MJ$ over $MJ$ we have that $X$ is weakly independent and $M = \Span(X) + MJ$. From the second equality we infer that $(M / \Span(X))J = M / \Span(X)$. Since the ideal $J$ is right $T$-nilpotent, we conclude that $M / \Span(X) = 0$, that is, $M = \Span(X)$. 
\end{prf}      

\begin{prp}\label{P01}Let $R$ be a ring.
\begin{enumerate} 
\item \cite[page 311]{NN} If all right $R$-modules are regularly weakly based, then $R$ is right perfect.
\item Let $J$ denote a Jacobson radical of $R$. If $R$ is right perfect, then all right $R$-modules are regularly weakly based if and only if all 
right modules over the semisimple ring $R / J$ are regularly weakly based.
\end{enumerate}
\end{prp}

\begin{prf} (1) follows readily from \Lmm{L04}, while (2) follows from \Lmm{L09}.
\end{prf}

\Prp{P01} reduces the characterization of rings over which all modules are regularly weakly based to a question whether all modules over a semisimple ring
are regularly weakly based. The answer to this seems  surprisingly non-trivial. 

We conclude the section with a straight consequence of \Prp{P01}. 

\begin{crl}\label{CrP01} Every module over a local perfect ring is regularly weakly based. 
\end{crl}

\section{Factoring out a finitely generating submodule}

It is easily seen that a module $M$ is weakly based if and only a factor $M/K$ is weakly based for a finitely generated submodule $K$ of $M$. The situation became less apparent when replacing weakly based with a regularly weakly based. We will apply this fact in the subsequent section. Before we proceed to its proof, we introduce the following notions (taken from \cite{HR2}).    

Let $M,N$ be modules, let $\phi: M \rightarrow N$ a homomorphism, and let $X$ be a subset of $M$. We say that $X$ \emph{lifts a subset $Y$ of $N$ \emph{via} $\phi$} provided that $\phi_{\restriction X}$ is a bijection onto $Y$. If $N$ is a quotient module of $M$, we say shortly that $X$ lifts $Y$, meaning that $X$ lifts $Y$ via the canonical projection.

Let $M$ be a module and let $X, Y$ be subsets of $M$. Let
\begin{equation*}
X^Y = \{x + \Span(Y) \mid x \in X\}.
\end{equation*} 
denote the image of the set $X$ in the canonical projection $M \to M / \Span(Y)$.

\begin{lmm} \label{T01}
Let $R$ be a ring, let $M$ be a right $R$-module and let $K$ be a finitely generated submodule of $M$. Then $M$ is regularly weakly based if and only if the factor module $M/K$ is regularly weakly based.
\end{lmm}

\begin{prf}
First suppose that the module $M$ is regularly weakly based. Let $\bar{X}$ be a generating set of $M/K$, and let $X$ be a subset of $M$ which lifts $\bar{X}$, i.e., $X^K = \bar X$. Then $X \cup K$ generates $M$, and since $M$ is regularly weakly based, $X \cup K$ contains a weak basis of $M$, say $Y$. Since $K$ is finitely generated, there is a finite subset $F$ of $Y$ such that $K \subseteq \Span(F)$. Put $Y_0 = Y \setminus F$. Since $Y$ is a weak basis of $M$, $Y_0^K$ is weakly independent in $M/K$. 

Since $Y$ generates $M$, the factor-module $M / (K + \Span(Y_0))$ is generated by $F^{K \cup Y_0}$. Since finitely generated modules are regularly weakly based, there is $F_0 \subseteq F$ that lifts a weak basis of $M / (K + \Span(Y_0))$. Since $Y_0$ is weakly independent over $F$, $K \subseteq \Span(F)$, and $F_0$ lifts a weak basis of $M / (K + \Span(Y_0))$, we conclude that $Y_0^K \cup F_0^K$ is a weak basis of $M/K$. Since $Y_0 \cup F_0 \subseteq X \cup K$, we infer that $Y_0 \cup F_0 \subseteq X$, whence  $Y_0^K \cup F_0^K \subseteq \bar X$. We have proved that the module $M/K$ is regularly weakly based. 
    
Now suppose that the factor-module $M / K$ is regularly weakly based. Let $X$ be a generating subset of $M$. Since $K$ is finitely generated, there is a finite subset $F$ of $X$ such that $K \subseteq \Span(F)$. The already proved implication gives that $M / \Span(F)$ is regularly weakly based. Thus we can pick a subset, $X_0$, of $X$ lifting a weak basis of $M / \Span(F)$. Observe that $F^{\Span(X_0)}$ generates the factor-module $M / \Span(X_0)$ and, since a finitely generated module is regularly weakly based, there is $F_0 \subseteq F$ lifting a weak basis of $M / \Span(X_0)$. We conclude that $X_0 \cup F_0$ is a weak basis of $M$ contained in $X$.           
\end{prf}

\section{Regularly weakly based modules over Dedekind domains} 

From now on we restrict ourselves to the case of Dedekind domains. Let $R$ be a Dedekind domain. We denote by $\Spec(R)$ the set of all non-zero prime (and thus, maximal) ideals of $R$. An $R$-module $T$ is $\emph{torsion}$ if $\Ann(m) \neq 0$ for any $m \in T$. Recall that any torsion $R$-module $T$ has a primary decomposition, that is, $T=\bigoplus_{\mathfrak{p} \in \Spec(R)}T_\mathfrak{p}$, where $T_\mathfrak{p}=\{m \in T \mid \Ann(m)=\mathfrak{p}^k \text{ for some $k$}\}$. We say that $T$ is $\mathfrak{p}$-primary if $T=T_\mathfrak{p}$. Alternatively, the $\mathfrak{p}$-primary part $T_\mathfrak{p}$ correspond naturally to the localization $T \otimes_{R} R_\mathfrak{p}$. In particular, we can view a $\mathfrak{p}$-primary $R$-module naturally as a module over the localization $R_\mathfrak{p}$. 

Let us recall a notion from abelian group theory which will prove useful in what follows. We say that a submodule $B$ of a $\mathfrak{p}$-primary module $T$ is \emph{basic} if $B$ is a pure submodule of $T$, $B$ is isomorphic to a direct sum of cyclic modules, and $T/B$ is divisible. As all these notions hold the same meaning independent of whether we view $T$ as an $R$-module or as an $R_\mathfrak{p}$-module, we can use \cite[Theorem 9.4]{KT} to infer that any $\mathfrak{p}$-primary module has a basic submodule (determined uniquely up to isomorphism).

Module $M$ is said to be \emph{bounded} if $IM=0$ for some non-zero ideal $I$.

\begin{lmm}\label{L05} Let $R$ be a Dedekind domain and let $T$ be a torsion $R$-module. If $T$ is regularly weakly based, then $T$ is bounded.
\end{lmm}
 
\begin{prf} Let $T$ be an unbounded torsion $R$-module. First suppose that $T$ is $\mathfrak{p}$-primary for some $\mathfrak{p} \in \Spec(R)$. We claim that there is a projection from $T$ onto a non-zero divisible module. In order to prove this, choose a basic submodule $B$ of $T$ (existence of which is discussed above). If $B \subsetneq T$, then $T/B$ is nonzero divisible and $T \rightarrow T/B$ is the desired projection. If $B=T$, then $T$ is a direct sum of $\mathfrak{p}$-primary cyclic modules of unbounded annihilators, and hence $T$ contains a submodule $S$ isomorphic to $\bigoplus_{n \in \NN}R/\mathfrak{p}^n$. It is well known that the indecomposable $\mathfrak{p}$-primary divisible $R$-module can be constructed as a direct limit of the system of inclusions $R/\mathfrak{p} \rightarrow R/\mathfrak{p}^2 \rightarrow R/\mathfrak{p}^3 \rightarrow \cdots$, and thus it is a quotient of $S$. As divisible $R$-modules are injective, this projection can be extended to the entire $T$. 

		We showed that there is a projection $\pi: T \rightarrow D$, where $D$ is non-zero divisible module. Denote by $K$ the kernel of $\pi$ and choose a generating set $X'$ of $D$. Since $D$ is divisible, there is a subset $X$ of $\mathfrak{p}T$ lifting $X'$ via $\pi$. Put $Z=X \cup K$ and note that $Z$ generates $T$. Suppose that $W \subseteq Z$ is a weak basis of $T$. By \cite[Corollary 3.3 and Lemma 5.2]{HR2}, any weak basis of $T$ lifts some basis of $T/\mathfrak{p}T$ over $\mathfrak{p}T$. Hence $W \subseteq K$, which is a contradiction to $W$ being a generating set.

		Let now $T$ be an unbounded (not necessarily $\mathfrak{p}$-primary) torsion $R$-module. Since regularly weakly based modules are closed under direct summands, the first part of this proof implies that $T_\mathfrak{p}$ is bounded for each $\mathfrak{p} \in \Spec(R)$. As $T$ is unbounded, there must be an infinite subset $\mathcal{P}$ of $\Spec(R)$ such that $T_\mathfrak{p} \neq 0$ for each $\mathfrak{p} \in \mathcal{P}$. If there is a non-zero divisible subgroup of $T$, it is a non-weakly based direct summand of $T$ (see \cite[Corollary 3.6]{HR2}). Thus $T$ is not regularly weakly based. We can thus assume that $T$ is reduced and apply \cite[Theorem 9]{K} to infer that there is a non-zero cyclic direct summand $C_\mathfrak{p}$ of $T_\mathfrak{p}$ for each $\mathfrak{p} \in \mathcal{P}$. Since $\mathcal{P}$ is infinite, we can pick a countable infinite sequence $\mathfrak{p}_n$, $n \in \go$ of pairwise distinct primes from $\eP$. It will suffice to show that $\bigoplus_{n \in \omega}C_{\mathfrak{p}_n}$ is not regularly weakly based. Fix a generator $x_n$ of $C_{\mathfrak{p}_n}$ and put $y_n = x_0 + x_1 + \cdots + x_n$ for each $n \in \go$. It follows easily that $\Span(y_m) \subsetneq \Span(y_n)$ whenever $m < n$, and so the generating set $\{y_n \mid n \in \omega\}$ of $\bigoplus_{n \in \omega}C_{\mathfrak{p}_n}$ does not contain a weak basis.
\end{prf}

\begin{lmm}\label{L13} Let $R$ be a Dedekind domain, and let $\mathfrak{p} \in \Spec(R)$. Every bounded $\mathfrak{p}$-primary $R$-module is regularly weakly based.
\end{lmm}

\begin{prf} Let $B$ be a bounded $\mathfrak{p}$-primary $R$-module. Then $B \mathfrak{p}^n = 0$ for some positive integer $n$ and $B$ can be naturally viewed as an $R / \mathfrak{p}^n$ module. Since $R$ is a Dedekind domain, the factor ring $R / \mathfrak{p}^n$ is local perfect, hence $B$ is regularly weakly based by \Crl{CrP01}. 
\end{prf}

Before proving the main lemma of the paper, we need the following auxiliary lemma:

\begin{lmm} \label{L06} Let $R$ be a Dedekind domain and let $N$ be a torsion-free $R$-module. If $N$ is an extension of a free module by a torsion bounded module, then $N$ is projective.
\end{lmm}
\begin{prf} Let $F$ be a free submodule of $N$ such that the factor-module $B = N/F$ is bounded torsion. Enumerate a free basis $X = \{ x_\ga \mid \ga < \gl \}$ of $F$ by an ordinal $\gl$ and put $F_\gb = \Span( \{ x_\ga \mid \ga < \gb \})$ for all $\gb < \gl$. For each $\ga < \gl$, let $N_\ga$ denote the smallest pure subgroup of $N$ containing $F_\ga$. It follows that $N = \bigcup_{\ga < \gl} N_\ga$ is a filtration of $N$ with $N_{\ga + 1} / N_\ga$ torsion-free for each $\ga < \gl$. Finitely generated torsion free modules over Dedekind domains are projective, see \cite[Theorem 6.3.23]{BK}, and therefore it will suffice to prove that all $N_{\ga+1} / N_\ga$ are finitely generated (and thus, projective). Indeed, then $N \simeq \bigoplus_{\ga < \gl} N_{\ga + 1} / N_\ga$ and so $N$ is projective.  

	Put $B_\alpha=N_\alpha/ F_\alpha$ for each $\alpha<\lambda$. As $ F_\ga = N_\alpha \cap F$ by the independence of $X$, we have the isomorphism $B_\ga = N_\ga / (N_\ga \cap F) \simeq (N_\ga + F) / F$ and so we can view naturally $B_\ga$ as a submodule of $B$. Denote by $Q$ the field of quotients of $R$. For each $\alpha<\lambda$, we obtain the following commutative diagram:
		$$
			\begin{CD}
				0 @>>> (N_\alpha+\Span(x_{\alpha}))/N_\alpha @>>> N_{\alpha+1}/N_\alpha @>>> B_{\alpha+1}/B_\alpha @>>> 0 \\
				& & @VV\simeq V  @VV\subseteq V @VV\subseteq V \\
				0 @>>> R @>>> Q @>>> Q/R @>>> 0
			\end{CD}
		$$
		Both exact sequences in the rows are given by the obvious quotient maps. For the maps in columns, the left most isomorphism follows from the fact, that $\Span(x_{\alpha}) \cap N_\alpha=0$, as $N_\alpha$ is the purification of $X_\alpha$, and $x_\alpha \not\in X_\alpha$. The middle inclusion is given by $N_{\alpha+1}/N_\alpha$ being torsion-free module of rank 1, and the right-most column map is induced by the two other ones. It is well known that $(Q/R)_\mathfrak{p}$ is uniserial for each $\mathfrak{p} \in \Spec(R)$. As $B_{\alpha+1}/B_\alpha$ is bounded, it has only finitely many non-zero $\mathfrak{p}$-primary parts, and as each of them is a bounded submodule of a uniserial module, they are all finitely generated. Therefore, $B_{\alpha+1}/B_\alpha$ is finitely generated. We conclude that $N_{\alpha+1}/N_\alpha$ is an extension of a cyclic module by a finitely generated module, hence
it is finitely generated. This finishes the proof.
\end{prf}

\begin{lmm}\label{L07} A regularly weakly based module over a Dedekind domain splits into a direct sum of a projective module and a bounded torsion module.
\end{lmm}
\begin{prf} Let $M$ be a regularly weakly based module over a Dedekind domain $R$. Let $T$ denote the torsion submodule of $M$, and let $F = M / T$ be the torsion-free quotient of $M$. If $F$ is projective, then $M$ decomposes as $T \oplus F$, and both the direct summands are regularly weakly based, in particular, 
$T$ is bounded torsion by Lemma~\ref{L05}.

		Suppose now that $F$ is not projective. Then we start with the following claim:
\begin{clm}\label{CL07} There is an ideal $\mathfrak{p} \in \Spec(R)$ and a subset $X$ of $M$ which lifts a basis of $M / M\mathfrak{p}$ over $M\mathfrak{p}$, such that $M / \Span(X)$ is not regularly weakly based.
\end{clm}		

\begin{cprf} We choose arbitrary $\mathfrak{p} \in \Spec(R)$ and a subset $X'$ of $T$ lifting a basis of $T/T\mathfrak{p}$. As $T$ is a pure submodule of $M$, we can extend $X'$ to a subset $X$ of $M$ containing $X'$ such that $X$ lifts a basis of $M/M \mathfrak{p}$. Put $Y=X \setminus X'$ and note that $Y^T$ lifts a basis of $F/F\mathfrak{p}$ over $F\mathfrak{p}$. By \cite[Lemma 7.1]{HR2}, $Y^T$ is a linearly independent subset of $F=M/T$, hence $\Span(Y^T)$ is free. Set $D = M/\Span(X)$. We claim that $D$ is not regularly weakly based. 
		
		If $D$ is torsion, then $M/(T+\Span(X)) \simeq F/\Span(Y^T)$ is torsion too. As $F$ is an extension of $\Span(Y^T)$ by $M/(T+\Span(X))$, the latter module is not bounded by \Lmm{L06}, otherwise $F$ would be projective. Hence, $D$ is also an unbounded torsion module, and, by \Lmm{L05}, $D$ is not regularly weakly based as desired. 
	
		Finally, suppose that $D$ is not torsion. In this case, choose any element $d \in D$ with $\Ann(d)=0$, and put $D'=D/d\mathfrak{p}$. Because $dR \simeq R$, we have that $d\mathfrak{p} \subsetneq dR \subseteq D$, and thus there is a submodule of $D'$ isomorphic to $R/\mathfrak{p}$, showing that the $\mathfrak{p}$-primary component of $D'$ is non-zero. Since $D=D\mathfrak{p}$, also $D'=D' \mathfrak{p}$. As the $\mathfrak{p}$-primary component of $D'$ is a pure submodule of $D'$, it is divisible by $\mathfrak{p}$, and therefore divisible. Hence, $D'$ contains a non-zero divisible direct summand, and thus $D'$ is not regularly weakly based by \cite[Corollary 3.6]{HR2}. As $d\mathfrak{p}$ is a finitely generated submodule of $D$, Lemma~\ref{T01} shows that $D$ is not regularly weakly based as desired. This concludes the proof of the claim.
\end{cprf}
		 
	We pick a generating set $Y'$ of $M/\Span(X)$ which does not contain any weak basis. As $M/\Span(X)$ is divisible by $\mathfrak{p}$, we can find a subset $Y$ of $\mathfrak{p}M$ lifting $Y'$ over $\Span(X)$. Then $X \cup Y$ is a generating set, which does not contain any weak basis of $M$. Indeed, any subset of $X \cup Y$ generating $M$ must contain the entire $X$ and $Y'$ does not contain any weak basis of $M/\Span(X)$.
\end{prf}

\begin{thm}\label{T02} Let $R$ be a Dedekind domain that is not a division ring. Then a regularly weakly based $R$-module splits into a direct sum of a finitely generated projective module and a bounded torsion module. 
\end{thm}

\begin{prf} Let $M$ be regularly weakly based module over a Dedekind domain $R$. Then $M = P \oplus B$, where $P$ is projective and $B$ a bounded torsion $R$-module, by \Lmm{L07}. Since $R$ is not a division ring, it is not perfect, indeed the only perfect domains are division rings. Applying \Lmm{L04} and the fact that regularly weakly based modules are closed under direct summands, we conclude that $P$ is finitely generated.   
\end{prf}

\begin{lmm}\label{L12} Let $R$ be a Dedekind domain, $F$ a finitely generated module, and $B$ a bounded $\mathfrak{p}$-primary module. Then $F \oplus B$ is regularly weakly based.
\end{lmm}

\begin{prf} Apply \Lmm{T01} and \Lmm{L13}.
\end{prf}

\begin{crl}\label{Cr1T02} Let $R$ be a discrete valuation ring and $M$ an $R$-module. Then $M$ is regularly weakly based if and only if $M \simeq F \oplus B$, where $F$ is finitely generated free module, and $B$ is bounded torsion module. 
\end{crl}

\begin{crl}\label{Cr2T02} Let $A$ an abelian group. If $A$ is regularly weakly based, then $A \simeq F \oplus B$, where $F$ is finitely generated free and $nB = 0$ for some positive integer $n$.
\end{crl}

\section{Closing remarks}

The remaining question is whether any bounded torsion module over a Dedekind domain is regularly weakly based. In other words, we ask whether all $R/I$-modules are regularly weakly based for any non-zero ideal $I$ of a Dedekind domain $R$. Since a non-zero ideals over Dedekind domains are products of prime ideals, $I = 
P_1^{n_1} \cdots P_k^{n_k}$, where $P_1,\dots,P_k$ are distinct prime ideals and $n_1,\dots,n_k$ are positive integers. The Jacobson radical of the ring $R / I$ corresponds to the ideal $(P_1 \cdots P_k) / I$ and it is clearly nilpotent. Applying \Lmm{L09} we can reduce the question to the case when $I$ is a product of distinct primes. In this case $R / I = R / (P_1 \cdots P_k) \simeq (R/ P_1) \times \cdots \times (R/P_k)$ is a product of fields, i.e, it is a commutative semisimple ring. Thus we arrived to a particular case of the question discussed at the end of Section~\ref{Perfect}. Let us formulate it as an open problem:

\begin{prb}\label{Prb1} Is every module over a semisimple ring regularly weakly based. In particular, is every module over a product of division rings (fields) regularly weakly based?    
\end{prb}

			The class of regularly weakly based modules is not closed under submodules in general. A counterexample can be obtained as follows. Let $R$ be a commutative Von Neumann regular ring with infinitely generated socle $S$ (e.g. an infinite product of fields). The regular module $R$, being finitely generated, is regularly weakly based. We show that the $R$-module $S$ is not. There is a submodule (and thus, a direct summand) $S'$ of $S$ of length $\aleph_0$, say $S' \simeq \bigoplus_{n \in \omega}S_n$, with $S_n$ simple for each $n \in \omega$. As $R$ is regular, $S_n$ has a direct complement $M_n$ in $R$ for each $n$. Choose a generator $x_n$ of $S_n$ for each $n \in \omega$ and put $y_n=x_0 + x_1 + \cdots + x_n$. We claim that $Y=\{y_n \mid n \in \omega\}$ is a generating set of $S'$ which does not contain any weak basis. As $M_0 \cap M_1 \cap \dots \cap M_{n-1} \not\subseteq M_n$, we conclude that $\Span(y_n) \subseteq \Span(y_m)$ for each $n \leq m$, and that $\Span(x_n) \subseteq \Span(y_n)$ for each $n \in \omega$. Hence, $Y$ generates $S'$, and as $S'$ is not finitely generated, $Y$ contains no weak bases of $S'$.

	\begin{prb}\label{Prb2}
		Is the class of regularly weakly based modules always closed under quotients?
	\end{prb}

\end{document}